\documentclass[lettersize,journal]{IEEEtran}
\usepackage{amsmath,amsfonts}
\usepackage{algorithmic}
\usepackage{algorithm}
\usepackage{array}
\usepackage[caption=false,font=normalsize,labelfont=sf,textfont=sf]{subfig}
\usepackage{textcomp}
\usepackage{stfloats}
\usepackage{url}
\usepackage{verbatim}
\usepackage{graphicx}
\usepackage{hyperref}
\usepackage{url}
\usepackage[most]{tcolorbox}
\definecolor{darkspringgreen}{rgb}{0.09, 0.45, 0.27}

\DeclareMathOperator*{\argmin}{arg\,min}
\usepackage[backend=biber, style=ieee]{biblatex}
\usepackage{url}
\usepackage{hyperref}
\begin{document}
\title{Stochastic Optimisation Framework using the Core Imaging Library and Synergistic Image Reconstruction Framework  for PET Reconstruction\footnote{https://agenda.infn.it/event/36860/contributions/230108/}}
\author{Evangelos Papoutsellis$^{1}$, Casper da Costa-Luis$^{2}$, Daniel Deidda$^{3}$, Claire Delplancke$^{4}$, Margaret Duff$^{2}$, Gemma Fardell$^{2}$, Ashley Gillman$^{5}$, Jakob S. J\o{}rgensen$^{6}$,  Zeljko Kereta$^{7}$, Evgueni Ovtchinnikov$^{2}$,  Edoardo Pasca$^{2}$,  Georg Schramm$^{8}$ and  Kris Thielemans$^{9}$ \\[10pt] 10th Conference on PET, SPECT, and MR Multimodal Technologies, Total Body and Fast Timing in Medical Imaging, 20-23 May 2024, Isola d'Elba, Italy 
\thanks{$^{1}$Finden Ltd, Rutherford Appleton Laboratory,  Harwell Campus, UK, $^{2}$Scientific Computing Department, Science \& Technology Facilities Council, Harwell Campus, UK, , $^{3}$ National Physical Laboratory, UK, $^{4}$\'Electricit\'e de France, Research and Development, $^{5}$
Australian e-Health Res. Ctr., CSIRO, Brisbane, Queensland, Australia, $^{6}$Department of Applied Mathematics and Computer Science, Technical University of Denmark,  $^{7}$Department of Computer Science, University College London, UK, $^{8}$Department of Imaging and Pathology, Division of Nuclear Medicine, KU Leuven, Leuven, Belgium, $^{9}$Institute of Nuclear Medicine, University College London, UK.

With thanks for discussions and contributions from: Matthias Ehrhardt, Tang Junqi, Laura Murgatroyd, Sam Porter, Imraj Singh and Robert Twyman. 

E. Pap acknowledges funding through the Innovate UK Analysis for Innovators (A4i) program ``Denoising of chemical imaging and tomography data (Project No. 10060435)". The development of CIL is supported by CCPi (EPSRC grant EP/T026677/1) and the Ada Lovelace Centre at STFC. The development of SIRF is funded by CCP SyneRBI (EPSRC grant EP/T026693/1). J.S.J. is supported by The Villum Foundation (Grant No. 25893) Z.K. is supported by the UK EPSRC grant EP/X010740/1. C.D. is supported by the ``PET++: Improving Localization, Diagnosis and Quantification in Clinical and Medical PET Imaging with Randomized Optimization’ EP/S026045/1."}}

\maketitle
\begin{abstract}
We introduce a stochastic framework into the open--source Core Imaging Library (CIL) which enables easy development of stochastic algorithms. Five such algorithms from the literature are developed, Stochastic Gradient Descent, Stochastic Average Gradient (-Am\'elior\'e), (Loopless) Stochastic Variance Reduced Gradient. We showcase the functionality of the framework with a comparative study against a deterministic algorithm on a simulated 2D PET dataset, with the use of the open-source Synergistic Image Reconstruction Framework. We observe that stochastic optimisation methods can converge in fewer passes of the data than a standard deterministic algorithm\footnote{https://agenda.infn.it/event/36860/contributions/230108/}. 
\end{abstract}

\begin{IEEEkeywords}
stochastic algorithms, positron emission tomography, image reconstruction, software and quantification 
\end{IEEEkeywords}

\section{Introduction}
\IEEEPARstart{I}{}terative reconstruction methods have been applied with
great success for solving challenging optimisation problems, such as total variation (TV) regularisation. Since iterative methods are computationally demanding due to the increasingly large data sizes, a range of stochastic optimisation algorithms have been proposed in the literature to reduce the computational effort.

In this work, we extend the optimisation functionality of the Core Imaging Library (CIL) \cite{jorgensen2021core, papoutsellis2021core} with a stochastic framework that enables developing a range of stochastic algorithms found in the literature: Stochastic Gradient Descent (SGD),  Stochastic Average Gradient (SAG) \cite{schmidt2017minimizing}, SAG- Am\'elior\'e (SAGA) \cite{defazio2014saga}, Stochastic Variance Reduced Gradient (SVRG) \cite{johnson2013accelerating} and Loopless SVRG (LSVRG) \cite{kovalev2020don}. 
These add to the Stochastic Primal-Dual Hybrid-Gradient (SPDHG) currently available in CIL, see \cite{papoutsellis2021core} and references therein. 
We demonstrate the use of this framework on a Positron Emission Tomography (PET) application, thanks to the combined use of the Synergistic Image Reconstruction Framework (SIRF) \cite{ovtchinnikov2020sirf}.

The developed framework allows for an easy comparison between stochastic gradient estimators. In this summary, we observe that stochastic optimisation methods can converge in fewer passes of the data than a deterministic benchmark.

\begin{figure}
    \centering
    \includegraphics[width=1\linewidth]{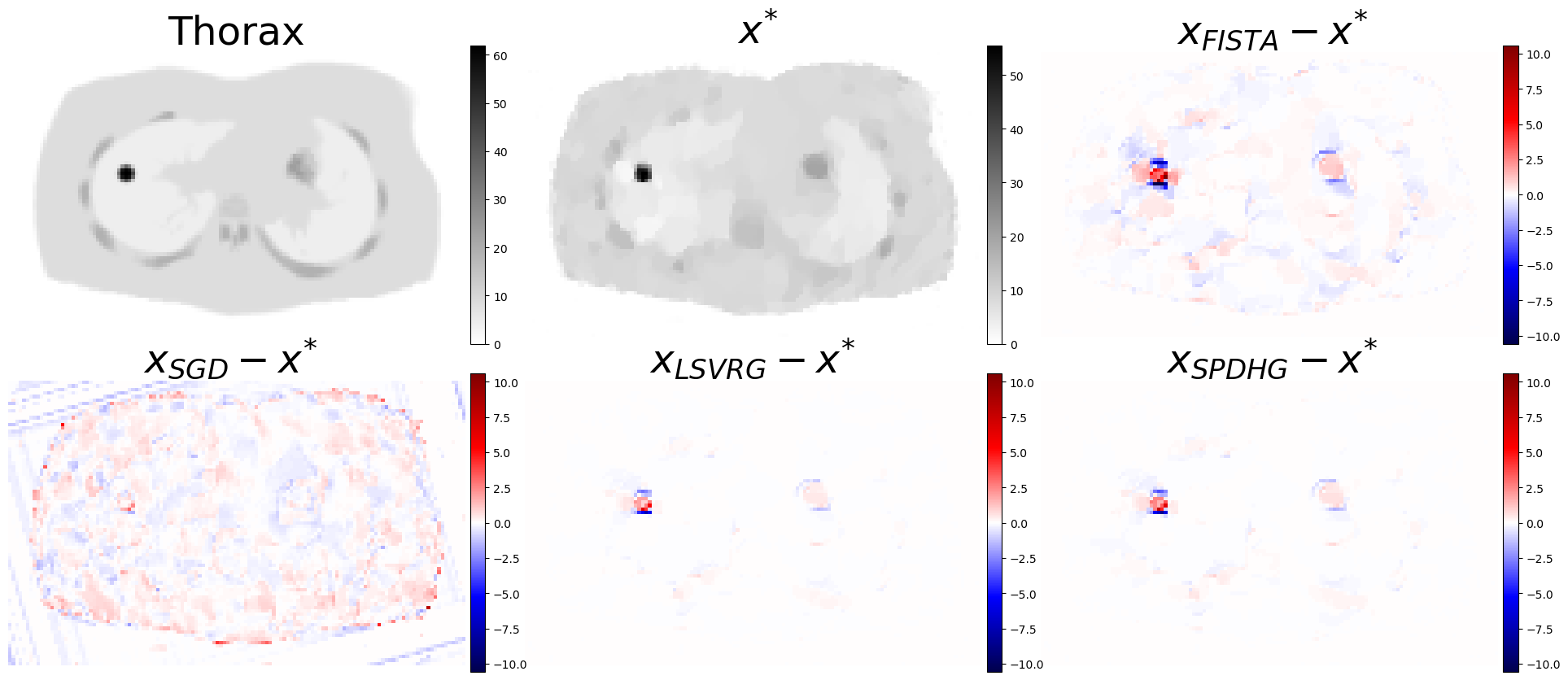}
    \caption{Simulated 2D $^{18}$F-FDG PET thorax dataset, reference solution $x^{*}$ and error plots for FISTA, Prox-SGD, Prox-LSVRG and SPDHG at 100 data passes.}
    \label{fig:data}
\end{figure}

\section{Stochastic Framework}

We consider optimisation problems of the form 
\begin{equation}\label{eq:optimisation_problem}
    x^*\!=\!\argmin_{x\in \mathbb{X}} \{F(x):= f(x)+g(x)\} \equiv\sum_{i=1}^{n}f_i(x)+g(x),
\end{equation}
where $\mathbb{X}$ is a finite dimensional space. Functions $f_i, f, g:\mathbb{X}\rightarrow \mathbb{R}$ for $i\in \{1,...,n\}$, are proper and convex, where $f_i$ are $L$-smooth and represent the fitness to the data. 
Regulariser $g$ has a proximal operator which either has a closed-form representation or can be efficiently solved up to some precision. 

Proximal gradient descent \cite{beck2009fast} (also known as ISTA, or forward-backward splitting) is a classical deterministic algorithm to solve \eqref{eq:optimisation_problem} by the iterations
\begin{equation}\label{eq:ISTA}
x_{k+1}=\text{prox}_{\gamma_k g} \left(x_k-\gamma_k\nabla f(x_k)\right), \quad k = 0, 1, 2, \dots
\end{equation}
for step size $\gamma_k$ and initial guess $x_0$. When $g\equiv 0$ this reduces to gradient descent.
Instead of computing the full gradient $\nabla f(x_k)$ in each iteration, stochastic optimisation algorithms employ an estimator $\tilde\nabla f(x_k)$, typically using the information of only one randomly selected function $f_i$. 

The stochastic framework in CIL consists of four components that can be combined in a plug-and-play fashion: i) functions providing stochastic estimators for the gradient of $f$, ii) \textit{sampling} methods which take in a set of probabilities $p_i$ for choosing each of the functions $f_i$, iii) a \textit{partitioner} to split up the data, defining the $f_i$'s, and, iv) algorithms to solve \eqref{eq:optimisation_problem}.

Thus far we have implemented 5 functions that provide stochastic estimators for the gradient of $f$.
When these are used in combination with GD or ISTA algorithms of CIL, they correspond to SGD (Prox-SGD), SAG (Prox-SAG) \cite{schmidt2017minimizing}, SAGA (Prox-SAGA) \cite{defazio2014saga}, SVRG (Prox-SVRG) \cite{johnson2013accelerating} and LSVRG (Prox-LSVRG)  \cite{kovalev2020don}. Due to our flexible design, the above stochastic estimators can also be combined with Nesterov-type accelerated algorithms, e.g., FISTA \cite{beck2009fast}, see \cite{Driggs2020} for more details.

\section{Methodology and Results}

For the numerical study, we use a simulated 2D $^{18}$F-FDG PET dataset from SIRF~\footnote{https://github.com/SyneRBI/SIRF\_data/tree/master/examples/PET/thorax\_single\_slice
}. Simulated Poisson noise is applied to the acquisition data. The data is partitioned into 32 subsets with equidistant projection views. Kullback-Leibler data fitting term is used for $ f_i$s and $f$. TV with a non-negativity constraint is used for the regulariser, $g=\alpha \text{TV}$, where $\alpha=0.1$. Results are shown in figure~\ref{fig:data}. 

\begin{figure}
    \centering
    \includegraphics[width=0.98\linewidth,trim=0 4mm 0 6mm, clip]{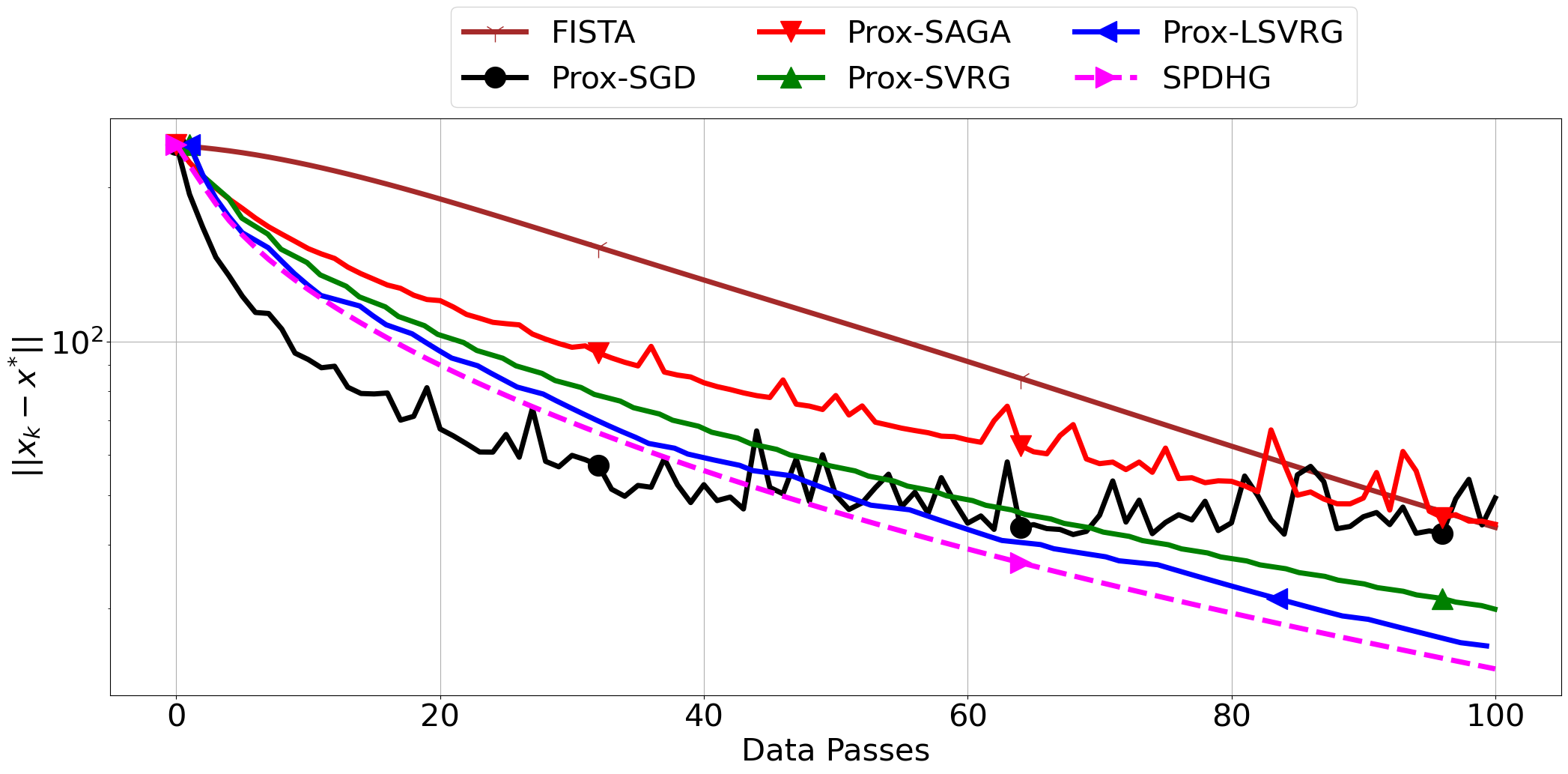}
    \caption{Distance from the optimal solution $x^*$ with respect to ``data passes".}
    \label{subopt-iterate}
\end{figure}

\begin{figure}
     \centering
     \includegraphics[width=0.98\linewidth,trim=0 4mm 0 6mm, clip]{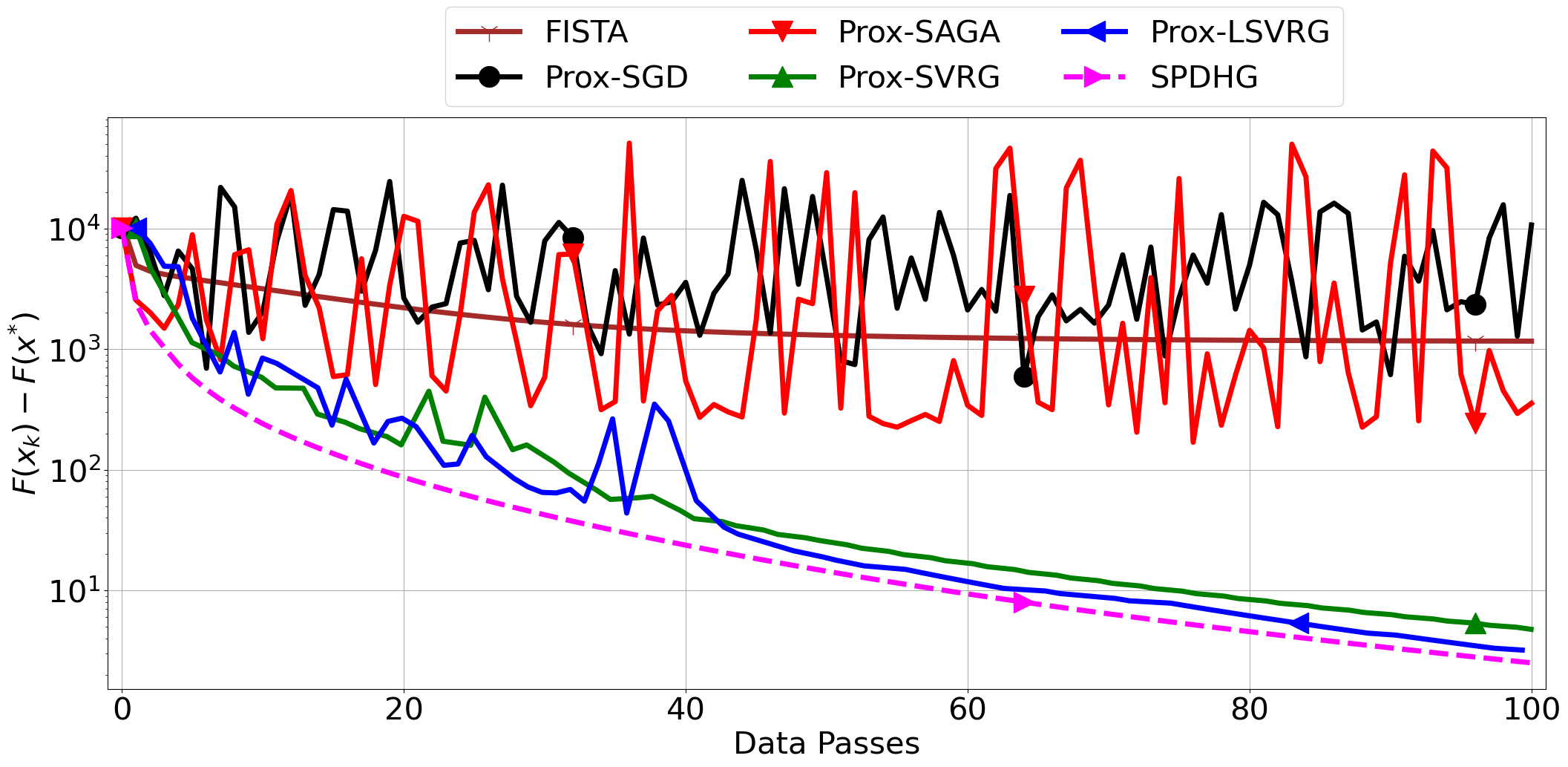}
     \caption{Distance from the optimal objective value $F(x^*)$  with respect to ``data passes".}
     \label{subopt-objectives}
 \end{figure}

 The optimal reconstruction $x^{*}$ was obtained using 500 data passes of SPDHG.
 All the algorithms are warm-started with one data pass of Prox-SGD. Functions $f_{i}$ are selected randomly with replacement. Algorithmic parameters such as step size, update frequency for Prox-SVRG, and probability for Prox-LSVRG were optimised using a parameter search. In figures \ref{subopt-iterate} and \ref{subopt-objectives}, we compare the different stochastic algorithms in CIL, and their performance with respect to ``data passes", i.e. how many times the algorithm has processed all the acquisition data in expectation.
 All the proposed stochastic algorithms have a faster convergence rate to the optimal solution than the deterministic FISTA.

\section{Discussion and Future Work}

This contribution describes an open-source framework that enables investigating a large variety of optimisation algorithms in many different contexts, including CT, MR, PET and SPECT image reconstruction. Presented results for PET show that stochastic algorithms in CIL can converge in fewer data passes than a deterministic counterpart. However, 
comparing the stochastic algorithms fairly will require a thorough investigation of algorithmic parameters, such as step size regimes. We leave this comparison for future work.

In addition, work is in progress to further empirically validate the stochastic framework by applying it to real PET data and to expand the versatility of the stochastic framework by extending its applicability to a wider array of stochastic algorithms, diverse imaging modalities, and integrating additional methodologies such as acceleration and pre-conditioning.

\printbibliography
\end{document}